
\input amssym.def
\input amssym
\input epsf
\magnification=1100
\baselineskip = 0.25truein
\lineskiplimit = 0.01truein
\lineskip = 0.01truein
\vsize = 8.5truein
\voffset = 0.2truein
\parskip = 0.10truein
\parindent = 0.3truein
\settabs 12 \columns
\hsize = 5.4truein
\hoffset = 0.4truein

\setbox\strutbox=\hbox{%
\vrule height .708\baselineskip
depth .292\baselineskip
width 0pt}
\font\caps=cmcsc10

\font\bigtenrm=cmr10 at 14pt

\def\sqr#1#2{{\vcenter{\vbox{\hrule height.#2pt
\hbox{\vrule width.#2pt height#1pt \kern#1pt
\vrule width.#2pt}
\hrule height.#2pt}}}}
\def\square{\mathchoice\sqr46\sqr46\sqr{3.1}6\sqr{2.3}4}

\centerline{\bigtenrm DETECTING LARGE GROUPS}

\tenrm
\vskip 14pt
\centerline{MARC LACKENBY}
\vskip 18pt

\tenrm
\centerline{\caps 1. Introduction}
\vskip 6pt

A group is known as {\sl large} if one of its finite index subgroups
has a free non-abelian quotient. Large groups have many interesting
properties, for example, super-exponential subgroup growth
and infinite virtual first Betti number. It is therefore
useful to be able to detect them in practice. In this paper,
we will show how one may deduce that a finitely presented group is
large using an array of different structures: its
profinite and pro-$p$ completions, its first $L_2$-Betti number
and the `homology growth' of its finite index subgroups.

The detection of large groups was the aim of [4], where
the author gave a characterisation of large
finitely presented groups in terms of the existence of
a sequence of finite index subgroups satisfying certain conditions.
In this paper, we start by deducing the following consequence.

\noindent {\bf Theorem 1.1.} {\sl Let $G$ and $K$
be finitely presented (discrete) groups that have isomorphic
profinite completions $\hat G$ and $\hat K$. Then $G$ is large if and only
if $K$ is large.}

In the above result, the term `isomorphic' can be
taken to mean `isomorphic as groups', since any
group isomorphism between profinite completions
$\hat G$ and $\hat K$ is automatically continuous.
We do not require that the isomorphism $\hat G \rightarrow \hat K$
be induced by a homomorphism $G \rightarrow K$.

One can also define a group to be {\sl $p$-large}, for some
prime $p$, if it contains a normal subgroup with index
a power of $p$ that has a free non-abelian quotient. In a
similar spirit to Theorem 1.1, we will prove the following.

\noindent {\bf Theorem 1.2.} {\sl Let $G$ and $K$
be finitely presented (discrete) groups that have isomorphic
pro-$p$ completions for some prime $p$. Then $G$ is
$p$-large if and only if $K$ is $p$-large.}

A sample application of Theorem 1.2 is to weakly
parafree groups, which are defined in terms of the
lower central series, as follows. Denote the $i$th term of the
lower central series of a group by $\gamma_i( \ )$.
A group is {\sl weakly parafree} if there is some non-trivial
free group $F$ with the `same' lower central series
as $G$. This means that there is an isomorphism 
$F/\gamma_i(F) \rightarrow G /\gamma_i(G)$,
for each positive integer $i$, and that these
isomorphisms are compatible with each other in the
obvious way. A group is known as {\sl parafree} if
it is weakly parafree and residually nilpotent.
Many interesting examples of parafree groups are given
in [1] and their properties are investigated in [2].
A consequence of Theorem 1.2 is the following result.

\noindent {\bf Theorem 1.3.} {\sl Let $G$ be a finitely
presented, weakly parafree group with $b_1(G) > 1$.
Then $G$ is large.}

Here, $b_1(G)$ denotes the first Betti number of $G$. This
has the following corollary.

\noindent {\bf Corollary 1.4.} {\sl Any finitely presented,
parafree group is either large or infinite cyclic.}

There are also applications of Theorem 1.2 to low-dimensional
topology, including the following.

\noindent {\bf Theorem 1.5.} {\sl If two closed 3-manifolds
$M_1$ and $M_2$ are topologically ${\Bbb Z}_p$-cobordant,
for some prime $p$, then $\pi_1(M_1)$ is $p$-large if and only
if $\pi_1(M_2)$ is $p$-large.}

Recall that two closed 3-manifolds 
$M_1$ and $M_2$ are {\sl topologically ${\Bbb Z}_p$-cobordant}
if there is a topological 4-manifold $X$ such that
$\partial X = M_1 \cup M_2$, and such that the inclusion of
each $M_i$ into $X$ induces isomorphisms of homology
groups with mod $p$ coefficients. Thus, Theorem 1.5 represents an unexpected link between two very different
areas of low-dimensional topology: the theory of finite covers
of 3-manifolds, and 4-dimensional topology. 
We will investigate this connection further in a forthcoming paper [9].

One of the goals of this paper is also to relate largeness
to $L_2$-Betti numbers. The first $L_2$-Betti number 
is defined in [13] for any finitely presented group $G$ and is denoted
here by $b_1^{(2)}(G)$.

\noindent {\bf Theorem 1.6.} {\sl Let $G$ be a finitely 
presented (discrete) group that is virtually residually $p$-finite
for some prime $p$, and such that
$b_1^{(2)}(G) > 0$. Then $G$ is large.}

This gives new examples of large groups. By applying
results of Shalom from [15], we obtain the following.

\noindent {\bf Corollary 1.7.} {\sl Let $G$ be a finitely
presented, non-amenable, discrete subgroup of ${\rm SO}(n,1)$ or
${\rm SU}(n,1)$, with $n \geq 2$ and with critical
exponent strictly less than 2. Then $G$ is large.}

This is a consequence of Theorem 1.6 because Shalom
showed that such a group $G$ has $b_1^{(2)}(G) > 0$
(Theorem 1.5 in [15]).
And since it is finitely generated and linear over a field
of characteristic zero, it is
virtually residually $p$-finite, for all but finitely many primes $p$
(Proposition 9 in Window 7 of [11]). Corollary 1.7 can be viewed as a generalisation
of the following important result of Cooper, Long and Reid (Theorem 1.3 of [3]) 
to Lie groups other than ${\rm SO}(3,1)$.

\noindent {\bf Theorem 1.8.} {\sl Let $G$ be a finitely generated,
discrete subgroup of ${\rm SO}(3,1)$ that is neither virtually
abelian nor cocompact. Then $G$ is large.}

Theorem 1.6 is proved using two results. The first
is due to L\"uck (Theorem 0.1 of [12]). It relates
$b_1^{(2)}(G)$ to the ordinary first Betti number $b_1(G_i)$
of finite index normal subgroups $G_i$.

\noindent {\bf Theorem 1.9.} (L\"uck [12]) {\sl Let $G$ be a finitely presented
group, and let $\{ G_i \}$ be a nested sequence of finite
index normal subgroups such that $\bigcap G_i = 1$. Then
$$\lim_{i \rightarrow \infty} {b_1(G_i) \over [G:G_i]}$$
exists and equals $b_1^{(2)}(G)$.}

This theorem is concerned with the growth rate of
$b_1(G_i)$ for finite index subgroups $G_i$. Recent
work of the author has instead focused on the growth rate
of homology with coefficients modulo some prime.
Let us fix some terminology. Let ${\Bbb F}_p$
be the field of order a prime $p$. For a group $G$,
let $d_p(G)$ be the dimension of the homology group $H_1(G; {\Bbb F}_p)$.

The second result forming the basis for Theorem 1.6 is 
the following, which is a consequence of the results of
the author in [4].

\noindent {\bf Theorem 1.10.} {\sl Let $G$ be a finitely presented
group with a surjective homomorphism $\phi \colon G \rightarrow
{\Bbb Z}$. Let $G_i = \phi^{-1}(i {\Bbb Z})$, and
let $p$ be a prime. Then 
\item{1.} $\lim_{i \rightarrow \infty} d_p(G_i) / [G:G_i]$ exists;
\item{2.} this limit is positive if and only if $d_p(G_i)$ is unbounded;
\item{3.} if the limit is positive, then $G$ is large.

}

Thus, fast growth of $d_p(G_i)$ as a function of the index
$[G:G_i]$ appears to be a strong
and useful property. We say that a nested sequence of finite index
subgroups $\{G_i \}$ of a group $G$ has {\sl linear growth of
mod $p$ homology} if $\inf_i d_p(G_i)/[G:G_i]$ is strictly positive.
A notable situation where this arises
is the following, which was the main result of [5].

\noindent {\bf Theorem 1.11.} {\sl Let $G$ be a lattice in ${\rm PSL}(2, {\Bbb C})$
satisfying one of the following:
\item{1.} $G$ contains a non-trivial torsion element, or
\item{2.} $G$ is arithmetic.

\noindent Then $G$ contains a strictly nested sequence of finite index
subgroups $\{ G_i \}$ with linear growth of mod $p$ homology,
for some prime $p$.

}

It seems very likely that these lattices are large. But it remains
unclear whether the conclusion of the theorem is strong enough to
imply this. However, the following theorem provides an affirmative
answer when $\{ G_i \}$ is the derived $p$-series for $G$. Recall
that this is a sequence of finite index subgroups $\{ D_i^{(p)}(G) \}$ defined recursively
by setting $D_0^{(p)}(G) = G$ and $D_{i+1}^{(p)}(G) =
[D_i^{(p)}(G),D_i^{(p)}(G)] (D_i^{(p)}(G))^p$.
Thus, $D_i^{(p)}(G)/D_{i+1}^{(p)}(G)$ is simply $H_1(D_i^{(p)}(G); {\Bbb F}_p)$.

\noindent {\bf Theorem 1.12.} {\sl Let $G$ be a finitely presented group,
and let $p$ be a prime. Suppose that the derived $p$-series for $G$ has linear
growth of mod $p$ homology. Then $G$ is $p$-large.}

This has implications for other series of finite index subgroups of $G$, for
example the lower central $p$-series, which is defined as follows.
The first term $\gamma^{(p)}_1(G)$ is $G$. The remaining terms are defined
recursively, setting $\gamma^{(p)}_{i+1}(G) = [\gamma^{(p)}_i(G), G](\gamma^{(p)}_i(G))^p$.

\noindent {\bf Corollary 1.13.} {\sl Let $G$ be a finitely presented group,
and let $p$ be a prime. Suppose that the lower central $p$-series for $G$ has linear
growth of mod $p$ homology. Then $G$ is $p$-large.}

The reason that 1.12 implies 1.13 is as follows. Each term
$D_i^{(p)}(G)$ of the derived $p$-series contains $\gamma_j^{(p)}(G)$
for all sufficiently large $j$. This is because the lower central
$p$-series of the finite $p$-group $G/D_i^{(p)}(G)$ terminates in the identity element,
since this is true for any finite $p$-group.
Moreover, we have the inequality
$${d_p(D_i^{(p)}(G)) - 1 \over [G:D_i^{(p)}(G)]} \geq
{d_p(\gamma_j^{(p)}(G)) - 1 \over [G:\gamma_j^{(p)}(G)]}.$$
This is an application of Lemma 3.3 in this paper, using the fact that
$\gamma_j^{(p)}(G)$ is normal in $D_i^{(p)}(G)$ and has index a power of $p$.
Thus, the assumption that the lower central $p$-series of $G$
has linear growth of mod $p$ homology implies that the same
is true of the derived $p$-series. Theorem 1.12 then implies
that $G$ is $p$-large, as required. It is clear from this
proof that versions of Corollary 1.13 apply to other series
of subgroups, for example, the dimension subgroups modulo ${\Bbb Z}_p$.

Theorem 1.12 is a consequence of a more general result, which we now describe.
An {\sl abelian $p$-series} for a group $G$ is a 
sequence of finite index subgroups $G = G_1 \triangleright G_2
\triangleright G_3 \triangleright \dots$ such that $G_i/G_{i+1}$ is
an elementary abelian $p$-group for each natural number
$i$. We investigate finitely presented groups $G$ having
an abelian $p$-series $\{ G_i \}$ which descends as fast as possible,
in the sense that the index $[G_i:G_{i+1}]$ is
(approximately) as big as it can be. Clearly,
the fastest possible descent occurs for the
derived $p$-series of a non-abelian free group $F$, of
rank $n$, say.
In this case,
$$[D_i^{(p)}(F):D_{i+1}^{(p)}(F)]
= p^{d_p(D_i^{(p)}(F))} = p^{[F:D_i^{(p)}(F)](n - 1)+1},$$
and so 
$$d_p(D_i^{(p)}(F)/D_{i+1}^{(p)}(F)) = {[F:D_i^{(p)}(F)](n - 1)+1}.$$
Thus, we say that an abelian $p$-series $\{ G_i \}$ has
{\sl rapid descent} if
$$\inf_i {d_p(G_i/G_{i+1}) \over [G:G_i]} > 0.$$

In Sections 4-7, we will prove the following theorems.

\noindent {\bf Theorem 1.14.} {\sl Let $G$ be a finitely presented group,
and let $p$ be a prime. Then the following are equivalent:
\item{1.} $G$ is large;
\item{2.} some finite index subgroup of $G$ 
has an abelian $p$-series with rapid descent.}

\vfill\eject
\noindent {\bf Theorem 1.15.} {\sl Let $G$ be a finitely presented group,
and let $p$ be a prime. Then the following are equivalent:
\item{1.} $G$ is $p$-large;
\item{2.} $G$ has an abelian $p$-series with rapid descent.}

These theorems represent a significant improvement upon the results in [4],
and can be viewed as the strongest theorems in this paper. They are interesting
for two reasons. Firstly, (2) in each theorem does not obviously imply that $G$
has a finite index subgroup with positive $b_1$, although this
is of course a consequence of the theorems. Secondly, the proof
of these results involves some genuinely new techniques. As in [4] and [6], topological
and geometric methods play a central role. But in this paper,
some basic ideas from the theory of error-correcting codes are
also used. In particular, we apply a generalisation of the so-called
`Plotkin bound' [14].

\noindent {\sl Acknowledgement.} I would like to thank Yehuda Shalom for pointing
out Corollary 1.7 to me.

\vskip 18pt
\centerline{\caps 2. Profinite completions and weakly parafree groups}
\vskip 6pt

Our goal in this section is to prove Theorems 1.1 - 1.5. Our starting
point is the following result, which is one of the main theorems in [4].

\noindent {\bf Theorem 2.1.} {\sl Let $G$ be a finitely
presented group. Then the following are equivalent:
\item{1.} $G$ is large;
\item{2.} there exists a sequence $G_1 \geq G_2 \geq \dots$
of finite index subgroups of $G$, each normal in $G_1$, such that
\itemitem{(i)} $G_i/G_{i+1}$ is abelian for all $i \geq 1$;
\itemitem{(ii)} $\lim_{i \rightarrow \infty} 
((\log [G_i : G_{i+1}]) / [G:G_i]) = \infty$;
\itemitem{(iii)} $\limsup_i (d(G_i/G_{i+1}) / [G:G_i])  > 0$.
}

Here, $d( \ )$ denotes the rank of a group, which is the minimal
size of a generating set.

\noindent {\bf Remark 2.2.} In the proof of (2) $\Rightarrow$ (1), 
one actually deduces that $G_i$ admits a surjective homomorphism
onto a non-abelian free group, for all sufficiently large $i$.
(See the comments after Theorem 1.2 in [4].)

Theorem 1.1 is a rapid consequence of the above result
and the following elementary facts, which follow immediately
from the definition of the profinite completion of a group.

\noindent {\bf Proposition 2.3.} {\sl Let $G$ and $K$ be finitely
generated (discrete) groups, and let $\phi \colon \hat G \rightarrow
\hat K$ be an isomorphism between their profinite completions. Then
the following hold.
\item{1.} There is an induced bijection (also denoted $\phi$) between the set of finite index
subgroups of $G$ and the set of finite index subgroups of $K$.
\item{2.} If $G_i$ is any finite index subgroup of $G$,
then $G_i$ is normal in $G$ if and only if $\phi(G_i)$ is
normal in $K$. In this case, there is an induced
isomorphism $G/G_i \rightarrow K/\phi(G_i)$, again denoted $\phi$.
\item{3.} If $G_i$ and $G_j$ are finite index subgroups
of $G$, then $G_i \subset G_j$ if and only if $\phi(G_i)
\subset \phi(G_j)$.
\item{4.} If $G_i \subset G_j$ are finite index subgroups
of $G$, and $G_i$ is normal in $G$, then
the isomorphism $\phi \colon G/G_i \rightarrow K/\phi(G_i)$
sends $G_j/G_i$ to $\phi(G_j)/\phi(G_i)$.
}

We can now prove the following.

\noindent {\bf Theorem 1.1.} {\sl Let $G$ and $K$
be finitely presented (discrete) groups that have isomorphic
profinite completions $\hat G$ and $\hat K$. Then $G$ is large if and only
if $K$ is large.}

\noindent {\sl Proof.} 
Let $\phi \colon \hat G \rightarrow \hat K$ be the given
isomorphism. 
Suppose that $G$ is large. It therefore contains a nested
sequence of finite index subgroups $G_i$ satisfying each
of the conditions in Theorem 2.1. These conditions are all
detectable by the profinite completion, as follows. 

Let $\tilde G_i$ be the intersection of
the conjugates of $G_i$. Proposition 2.3 (1) gives finite
index subgroups $\phi(G_i)$ and $\phi(\tilde G_i)$ of $K$,
which we denote by $K_i$ and $\tilde K_i$ respectively.
By 2.3 (2), $\tilde K_i$ is normal in $K$. By 2.3 (3),
$\tilde K_i$ is contained in $K_i$. By 2.3 (2) and 2.3 (4), there
is an isomorphism $G/\tilde G_i
\rightarrow K/\tilde K_i$ which takes 
$G_j/\tilde G_i$ to $K_j/\tilde K_i$ for
any $j \leq i$. The normality of $G_i$ in $G_1$
is equivalent to the normality of $G_i/\tilde G_i$
in $G_1/\tilde G_i$. Thus, $K_i$ is normal in $K_1$.
The isomorphism $G/\tilde G_{i+1} \rightarrow
K/\tilde K_{i+1}$ takes $G_i/\tilde G_{i+1}$
and $G_{i+1}/\tilde G_{i+1}$ onto $K_i/\tilde K_{i+1}$
and $K_{i+1}/\tilde K_{i+1}$ respectively. Hence,
$K_i/K_{i+1}$ is isomorphic to $G_i/G_{i+1}$.
Thus, the sequence $K_i$ satisfies the conditions
of Theorem 2.1. So, $K$ is large. $\square$

\noindent {\bf Remark 2.4.} A modified version of Proposition 2.3
holds, where the phrase `profinite completions'
is replaced by `pro-$p$ completions for some prime $p$',
and `finite index subgroup(s)' is
replaced throughout by `subnormal subgroup(s) with index
a power of $p$'.

We now embark upon the proof of Theorem 1.2. For this, we need
the following variant of Theorem 2.1.

\noindent {\bf Theorem 2.5.} {\sl Let $G$ be a finitely
presented group and let $p$ be a prime. Then the following are equivalent:
\item{1.} $G$ is $p$-large;
\item{2.} there exists a sequence $G_1 \geq G_2 \geq \dots$
of subgroups of $G$, each with index a power of $p$ in $G$,
such that $G_1$ is normal in $G$, and each $G_i$ is normal in $G_1$, and where
the following hold:
\itemitem{(i)} $G_i/G_{i+1}$ is abelian for all $i \geq 1$;
\itemitem{(ii)} $\lim_{i \rightarrow \infty} 
((\log [G_i : G_{i+1}]) / [G:G_i]) = \infty$;
\itemitem{(iii)} $\limsup_i (d(G_i/G_{i+1}) / [G:G_i])  > 0$.
}

\noindent {\sl Proof.} (1) $\Rightarrow$ (2): Since $G$ is $p$-large,
some finite index normal subgroup $G_1$, with index a power of $p$,
admits a surjective homomorphism $\phi$
onto a non-abelian free group $F$. Define the following subgroups
of $F$ recursively. Set $F_1 = F$, and let $F_{i+1} = 
[F_i,F_i](F_i)^{p^i}$. Set $G_i = \phi^{-1}(F_i)$.
Then it is trivial to check that the conditions 
of (2) hold. 

(2) $\Rightarrow$ (1): By Theorem 2.1 and Remark 2.2, some $G_i$
admits a surjective homomorphism $\phi$ onto a non-abelian free group $F$. 
By assumption, $G_i$ is subnormal in $G$ and has index a power of
$p$. Set $\tilde G_i$ to be the intersection of the conjugates of $G_i$.
Then $\tilde G_i$ is normal in $G$ and also has index a power of $p$.
The restriction of $\phi$ to $\tilde G_i$ is a surjective homomorphism
onto a finite index subgroup of $F$, which is therefore free non-abelian.
Thus, $G$ is $p$-large.
$\square$

\noindent {\bf Theorem 1.2.} {\sl Let $G$ and $K$
be finitely presented (discrete) groups that have isomorphic
pro-$p$ completions for some prime $p$. Then $G$ is
$p$-large if and only if $K$ is $p$-large.}

\noindent {\sl Proof.} This is very similar to the proof
of Theorem 1.1. Suppose that $G$ is $p$-large. It
therefore has a sequence of subgroups $G_i$ where the
conclusions of Theorem 2.5 hold. Using Remark 2.4,
$K$ also has such a sequence of subgroups. Thus, by
Theorem 2.5, $K$ is $p$-large. $\square$

Our aim now is to prove Theorem 1.3.

\noindent {\bf Theorem 1.3.} {\sl Let $G$ be a finitely
presented, weakly parafree group with $b_1(G) > 1$.
Then $G$ is large.}

As stated in the Introduction, this has the following
corollary.

\noindent {\bf Corollary 1.4.} {\sl Any finitely presented,
parafree group is either large or infinite cyclic.}

\noindent {\sl Proof.} Any parafree group $G$ with $b_1(G) \leq 1$
is infinite cyclic. $\square$

Theorem 1.3 is a consequence of a stronger result concerning
weakly $p$-parafree groups, which we now define.
A group $G$ is {\sl weakly $p$-parafree}, for some prime $p$, if there is some
non-trivial free group $F$ such that $G/\gamma_i^{(p)}(G)$ is
isomorphic to $F/\gamma_i^{(p)}(F)$ for each $i \geq 1$.
Recall that $\gamma_i^{(p)}( \ )$ denotes the lower central
$p$-series of a group.

\noindent {\bf Proposition 2.6.} {\sl A weakly parafree
group is weakly $p$-parafree for each prime $p$.}

\noindent {\sl Proof.} By assumption, there is an isomorphism between 
$G/\gamma_i(G)$ and $F/\gamma_i(F)$ for some non-trivial free group $F$. This
induces an isomorphism between the lower central $p$-series
of $G/\gamma_i(G)$ and $F/\gamma_i(F)$. Hence,
$${G/\gamma_i(G) \over \gamma_i^{(p)}(G/\gamma_i(G))}
\cong {F/\gamma_i(F) \over \gamma_i^{(p)}(F/\gamma_i(F))}.$$
But $\gamma_i^{(p)}(G/\gamma_i(G))$ is isomorphic to
$\gamma_i^{(p)}(G)/\gamma_i(G)$, since $\gamma_i(G)$ is
contained in $\gamma_i^{(p)}(G)$. So, the left-hand side is
isomorphic to $G/\gamma^{(p)}_i(G)$. Similarly, the right-hand
side is isomorphic to $F/\gamma^{(p)}_i(F)$. Thus, we obtain an
isomorphism between $G/\gamma^{(p)}_i(G)$ and
$F/\gamma^{(p)}_i(F)$. $\square$

Note that, when defining the weakly $p$-parafree group $G$,
we do not make the assumption that the isomorphisms 
$G /\gamma_i^{(p)}(G) \rightarrow F/\gamma_i^{(p)}(F)$ are
{\sl compatible}. This means that the following
diagram commutes, for each $i \geq 2$:
$$

\matrix{
G / \gamma_i^{(p)}(G) & \longrightarrow & F / \gamma_i^{(p)}(F) \cr
\Big\downarrow && \Big\downarrow \cr
G / \gamma_{i-1}^{(p)}(G) & \longrightarrow & F / \gamma_{i-1}^{(p)}(F). \cr}
$$
Here, the horizontal arrows are the given isomorphisms
and the vertical maps are the obvious quotient homomorphisms.
However, we can assume this, with no loss, as the following
lemma implies.

\noindent {\bf Lemma 2.7.} {\sl Let $G$ and $K$ be finitely
generated groups. Suppose that, for each integer $i \geq 1$, there is an
isomorphism $\theta_i \colon G/\gamma_i^{(p)}(G) \rightarrow K/\gamma_i^{(p)}(K)$.
Then there is a collection of such isomorphisms that
are compatible.}

\noindent {\sl Proof.} For each $j \leq i$, $\theta_i$ restricts to an isomorphism between
$\gamma_j^{(p)}(G)/\gamma_i^{(p)}(G)$ and $\gamma_j^{(p)}(K)/\gamma_i^{(p)}(K)$, and hence, 
quotienting $G/\gamma_i^{(p)}(G)$ by $\gamma_j^{(p)}(G)/\gamma_i^{(p)}(G)$,
we obtain an isomorphism $\theta_{i,j} \colon
G/\gamma^{(p)}_j(G) \rightarrow K/\gamma^{(p)}_j(K)$.
As $G/\gamma^{(p)}_2(G)$ is finite, some $\theta_{i,2}$
occurs infinitely often. Take this to be the given isomorphism
$\phi_2 \colon G / \gamma^{(p)}_2(G) \rightarrow K / \gamma^{(p)}_2(K)$, and only consider
those $\theta_i$ for which $\theta_{i,2} = \phi_2$. Among these,
some $\theta_{i,3}$ occurs infinitely often. Define this to be
$\phi_3$, and so on. Then the $\phi_i$ form the required 
compatible collection of isomorphisms. $\square$

The above lemma is elementary and well-known, as is the following
result. They are included for the sake of completeness.

\noindent {\bf Lemma 2.8.} {\sl Let $G$ and $K$ be finitely
generated groups, and let $p$ be a prime. Then the following are equivalent:
\item{1.} the pro-$p$ completions of $G$ and $K$ are isomorphic;
\item{2.} for each $i \geq 1$, there is an
isomorphism $G/\gamma_i^{(p)}(G) \rightarrow K/\gamma_i^{(p)}(K)$.

}

\noindent {\sl Proof.} (1) $\Rightarrow$ (2): 
Let $\hat G_{(p)}$ denote the pro-$p$ completion of $G$.
Then, for each $i \geq 1$, $\hat G_{(p)}/\gamma_i^{(p)}(\hat G_{(p)})$
is isomorphic to $G/\gamma_i^{(p)}(G)$. The claim follows immediately.

(2) $\Rightarrow$ (1): The pro-$p$ completion $\hat G_{(p)}$ can be expressed
as the inverse limit of $\dots \rightarrow G/\gamma_2^{(p)}(G) \rightarrow G/\gamma_1^{(p)}(G)$. 
This is because the kernel of any homomorphism of $G$ onto a finite $p$-group
contains $\gamma_i^{(p)}(G)$ for all sufficiently large $i$. Suppose now
that, for each $i \geq 1$, there is an
isomorphism $G/\gamma_i^{(p)}(G) \rightarrow K/\gamma_i^{(p)}(K)$.
According to Lemma 2.7, these isomorphisms can be
chosen compatibly. This implies there is an isomorphism between the
inverse limits:
$$\lim_{\leftarrow} G/\gamma_i^{(p)}(G) \cong \lim_{\leftarrow} K/\gamma_i^{(p)}(K).$$
Thus, $\hat G_{(p)}$ and $\hat K_{(p)}$ are isomorphic. $\square$

Setting $K$ to be a free group in Lemma 2.8 gives
the following characterisation of
weakly $p$-parafree groups in terms of pro-$p$ completions.

\noindent {\bf Corollary 2.9.} {\sl Let $G$ be a finitely
generated (discrete) group and let $p$ be a prime. Then
$G$ is weakly $p$-parafree if and only if its pro-$p$
completion is isomorphic to the pro-$p$ completion of
a non-trivial free group.}

Thus, Theorem 1.3 is a consequence of the following.

\noindent {\bf Theorem 2.10.} {\sl Let $G$ be a finitely
presented group that is weakly $p$-parafree for some prime $p$.
Suppose that $d_p(G) > 1$. Then $G$ is $p$-large.}

\noindent {\sl Proof.} The assumption that $G$
is weakly $p$-parafree implies that the pro-$p$
completion of $G$ is isomorphic to the pro-$p$
completion of a free group $F$, by Corollary 2.9. Since $d_p(G) >1$,
$F$ is a non-abelian free group. In particular,
it is $p$-large. Thus, by Theorem 1.2, $G$ is $p$-large.
$\square$

We close this section with a topological application of
Theorem 1.2.

\noindent {\bf Theorem 1.5.} {\sl If two closed 3-manifolds
$M_1$ and $M_2$ are topologically ${\Bbb Z}_p$-cobordant,
for some prime $p$, then $\pi_1(M_1)$ is $p$-large if and only
if $\pi_1(M_2)$ is $p$-large.}

This is an immediate consequence of Theorem 1.2 and the following
result.

\noindent {\bf Theorem 2.11.} {\sl If two closed 3-manifolds
are topologically ${\Bbb Z}_p$-cobordant, then the pro-$p$
completions of their fundamental groups are isomorphic.}

We will prove this theorem in a forthcoming paper [9], where we will
develop further connections between 3-manifolds, 4-manifolds,
and the pro-$p$ completions of their fundamental groups.

\vfill\eject
\centerline{\caps 3. Homology growth in cyclic covers}
\vskip 6pt

The goal of this section is to prove the following result
and then Theorem 1.6.

\noindent {\bf Theorem 1.10.} {\sl Let $G$ be a finitely presented
group with a surjective homomorphism $\phi \colon G \rightarrow
{\Bbb Z}$. Let $G_i = \phi^{-1}(i {\Bbb Z})$, and
let $p$ be a prime. Then 
\item{1.} $\lim_{i \rightarrow \infty} d_p(G_i) / [G:G_i]$ exists;
\item{2.} this limit is positive if and only if $d_p(G_i)$ is unbounded;
\item{3.} if the limit is positive, then $G$ is large.

}

We will need the following lemma.

\noindent {\bf Lemma 3.1.} {\sl Let $k$ be a non-negative integer
and let $f \colon {\Bbb N}_{>0} \rightarrow {\Bbb R}$
be a function satisfying
$$|f(i+j) - f(i) - f(j)| \leq k,\eqno{(\ast)}$$
for any $i, j \in {\Bbb N}$. Then
\item{1.} $\lim_{i \rightarrow \infty} f_i / i$ exists;
\item{2.} this limit is non-zero if and only if $f_i$ is unbounded.

}

\noindent {\sl Proof.}
Note first that a trivial induction establishes that $f(i) \leq i(f(1)+ k)$
for each positive $i \in {\Bbb N}$.
Hence $\limsup_i f(i)/i$ is finite, $M$ say.

\noindent {\sl Claim.} Suppose that $f(m) > 2k$, for some positive $m \in {\Bbb N}$. 
Then $\liminf_i f(i)/i > 0$.

We will show that
$f(nm) > (n+1)k$, for each positive $n \in {\Bbb N}$,
by induction on $n$. The induction starts trivially.
For the inductive step, note that
$f((n+1)m) \geq f(m) + f(nm) - k > 2k + (n+1)k - k = (n+2)k$.
This establishes the inequality. The claim now follows by noting
that if $i = nm+r$, for $0 \leq r < m$, then
$$|f(i) - f(nm)| \leq k+ \max_{1 \leq r < m} |f(r)|.$$

Now let $g(i) = Mi - f(i)$. By the definition of $M$, 
$\liminf_i g(i)/i = 0$. Now, $g$ satisfies $(\ast)$ and so
by the claim, $g$ is bounded above by $2k$. Hence,
$f(i) \geq Mi - 2k$. Thus, $\liminf_i f(i)/i = M$.
This proves (1). To prove (2), note that if $f$ satisfies
$(\ast)$, then so does $-f$. Thus, applying the claim, we
deduce that either $|f(i)| \leq 2k$ for all positive $i$ or
$\lim_{i \rightarrow \infty} f(i)/i$ is non-zero.
$\square$

\noindent {\sl Proof of Theorem 1.10.} We claim that there is a non-negative integer $k$
such that, for all $i, j \geq 1$,
$$|d_p(G_{i+j}) - d_p(G_i) - d_p(G_j)| \leq k.$$
This and Lemma 3.1 will then imply (1) and (2).

Let $K$ be a finite connected 2-complex with fundamental group $G$.
We may find a map $F \colon K \rightarrow S^1$ so that
$F_\ast \colon \pi_1(K) \rightarrow \pi_1(S^1)$ is $\phi \colon G \rightarrow {\Bbb Z}$.
Let $b$ be a point in $S^1$. Then, after a small homotopy,
we may assume that $F^{-1}(b)$ is a finite graph $\Gamma$, that
a regular neighbourhood of $\Gamma$ is a copy of $\Gamma \times I$
and that the restriction of $F$ to this neighbourhood is
projection onto the $I$ factor, followed by inclusion of $I$ into $S^1$. Let $\overline K$ be the
result of cutting $K$ along $\Gamma$. Let $K_i$ be 
the $i$-fold cover of $K$ corresponding to $G_i$.
Then $K_i$ can be obtained from $i$ copies of $\overline K$
glued together in a circular fashion. Cut $K_i$ along
one of the copies of $\Gamma$ in $K_i$, and let $\overline K_i$
be the result. 
The Mayer-Vietoris sequence (applied to the decomposition of $K_i$
into $\overline K_i$ and $\Gamma \times I$) gives the following inequalities:
$$- d_p(\Gamma) \leq d_p(K_i) - d_p(\overline K_i) \leq |\Gamma|.$$
Similarly, since the disjoint union of $\overline K_i$
and $\overline K_j$ is obtained by cutting $\overline K_{i+j}$ along a copy of $\Gamma$, we have
$$-d_p(\Gamma) \leq d_p(\overline K_{i+j}) - d_p(\overline K_i) - d_p(\overline K_j) \leq
|\Gamma|.$$
Since $d_p(K_i) = d_p(G_i)$, the claim now follows, letting
$k = 2 d_p(\Gamma) + 2 |\Gamma|$.

Conclusion (3) is a consequence of the following result 
(Theorem 1.2 of [4]), setting $H_i = G$, $J_i = G_i$ and
$K_i = [G_i,G_i](G_i)^p$. $\square$

\noindent {\bf Theorem 3.2.} {\sl Let $G$ be a finitely presented group
and suppose that, for each natural number $i$, there is a triple
$H_i \geq J_i \geq K_i$ of finite index normal subgroups such that
\item{(i)} $H_i / J_i$ is abelian for all $i$;
\item{(ii)} $\lim_{i \rightarrow \infty} ((\log [H_i:J_i])/[G:H_i]) = \infty$;
\item{(iii)} $\inf_i (d(J_i/K_i)/[G:J_i]) > 0$.

\noindent Then, $K_i$ admits a surjective homomorphism onto a non-abelian free group,
for all sufficiently large $i$.}

To prove Theorem 1.6, we need one more fact,
which is well known. It appears as Proposition 3.7 in [5],
for example.

\noindent {\bf Lemma 3.3.} {\sl Let $G$ be a finitely
generated group and let $K$ be a normal subgroup with
index a power of a prime $p$. Then
$$d_p(K) \leq (d_p(G) - 1) [G:K] + 1.$$
}

We can now prove Theorem 1.6.

\noindent {\bf Theorem 1.6.} {\sl Let $G$ be a finitely presented 
(discrete) group, that is virtually residually $p$-finite, for some
prime $p$, and such that $b_1^{(2)}(G) > 0$. Then $G$ is large.}

\noindent {\sl Proof.} Since $G$ is virtually residually $p$-finite,
it has a finite index normal subgroup $G_1$ that is residually
$p$-finite. Thus, $G_1$ contains a nested sequence of normal subgroups $G_i$,
each with index a power of $p$, such that
$\bigcap_i G_i = 1$. By the multiplicativity of $b_1^{(2)}$
(Theorem 1.7 (1) of [13]), $b_1^{(2)}(G_1) = b_1^{(2)}(G) [G:G_1] > 0$.
By L\"uck's theorem (Theorem 1.9),
$\lim_{i \rightarrow \infty} b_1(G_i)/[G_1:G_i]$ exists
and equals $b_1^{(2)}(G_1)$, which is positive. Hence,
by relabelling the $G_i$, we may assume that
$b_1(G_1) > 0$. Let $\phi \colon G_1 \rightarrow {\Bbb Z}$ be
a surjective homomorphism,
and let $K_i = \phi^{-1}(p^i {\Bbb Z})$.

We claim that $\liminf_i d_p(K_i)/[G:K_i]$
is positive. Consider the subgroups $G_i \cap K_i$. Each
is a finite index normal subgroup of $G_1$ and their intersection is
the identity.
Hence, by Theorem 1.9, $\lim_{i \rightarrow \infty}
b_1(G_i \cap K_i) / [G_1: G_i \cap K_i]$ exists and is
positive. In particular, $\liminf_i d_p(G_i \cap K_i) / [G: G_i \cap K_i]$
is positive. Now, the quotient $K_i / (G_i \cap K_i)$ is
isomorphic to $K_i G_i / G_i$, which is a subgroup of $G_1/G_i$,
and so has order a power of $p$. Hence, by Lemma 3.3,
$${d_p(G_i \cap K_i) \over [G:G_i \cap K_i]}
\leq {(d_p(K_i) - 1)[K_i:G_i \cap K_i] + 1 \over [G:G_i \cap K_i]}.$$
Therefore, $\liminf_i d_p(K_i)/[G:K_i]$ is positive, as claimed.

In particular, $d_p(K_i)$ is unbounded. Thus, by Theorem 1.10,
$G$ is large. $\square$

\vfill\eject
\centerline{\caps 4. Cocycle size and Property $(\tau)$}
\vskip 6pt

Most of the remainder of the paper is devoted to the proof of
Theorems 1.14 and 1.15.

\noindent {\bf Theorem 1.14.} {\sl Let $G$ be a finitely presented group
and let $p$ be a prime. Then the following are equivalent:
\item{1.} $G$ is large;
\item{2.} some finite index subgroup of $G$ 
has an abelian $p$-series with rapid descent.}

\noindent {\bf Theorem 1.15.} {\sl Let $G$ be a finitely presented group,
and let $p$ be a prime. Then the following are equivalent:
\item{1.} $G$ is $p$-large;
\item{2.} $G$ has an abelian $p$-series with rapid descent.}

The difficult direction in each of these theorems is $(2) \Rightarrow (1)$.
For this, one needs a method for proving that a finitely
presented group is large or $p$-large. Various techniques have been
developed with this aim. The one we will use deals with
the `relative size' of cocycles on 2-complexes.

Let $K$ be a connected finite 2-complex with fundamental
group $G$. Let $K_i$ be the covering
space corresponding to a finite index subgroup $G_i$.
The key to our approach is to consider cellular 1-cocycles
on $K_i$ representing non-trivial elements
of $H^1(K_i; {\Bbb F}_p)$.

For a cellular
1-dimensional cocycle $c$ on $K_i$, let its {\sl support} ${\rm supp}(c)$
be those 1-cells with non-zero evaluation under $c$.
For an element $\alpha \in H^1(K_i; {\Bbb F}_p)$,
consider the following quantity, which was defined in [6]. 
The {\sl relative size} of $\alpha$ is
$${\rm relsize}(\alpha) =
{\min \{ |{\rm supp}(c)|: c \hbox{ is a cellular cocycle 
representing } \alpha \} \over 
\hbox{Number of 1-cells of } K_i}.$$

The following result was proved in [6] and is
central to our approach.

\noindent {\bf Theorem 4.1.} {\sl Let $K$ be a finite connected
2-complex, and let $\{ K_i \rightarrow K\}$ be a
collection of finite-sheeted covering spaces.
Suppose that $\{ \pi_1(K_i) \}$ has linear
growth of mod $p$ homology for some prime $p$.
Then one of the following must hold:
\item{(i)} $\pi_1(K_i)$ is $p$-large for infinitely many $i$, or
\item{(ii)} there is some $\epsilon > 0$ such that
the relative size of any non-trivial class in $H^1(K_i; {\Bbb F}_p)$
is at least $\epsilon$, for all $i$.

}

This is a slightly modified version of Theorem 6.1 of [6].
In that result, it is not explicitly stated that
$\pi_1(K_i)$ is $p$-large for infinitely many $i$,
merely that $\pi_1(K)$ is large. But the proof does indeed
give a normal subgroup of $\pi_1(K_i)$ with index a
power of $p$ that has a free non-abelian quotient.

In this section, we will relate the relative size of
cocycles to Property $(\tau)$. While not directly needed
in the remainder of the paper, this is a potentially
important link.

We now recall the definition of Property $(\tau)$.
Let $G$ be a group with a finite generating set $S$.
Let $\{G_i \}$ be a collection of finite index subgroups
of $G$. Let $X_i = X(G/G_i; S)$ be the Schreier coset
graph of $G/G_i$ with respect to the generating set $S$.

The {\sl Cheeger constant} $h(X_i)$ is defined to be
$$h(X_i) = \min \left \{ {|\partial A| \over |A|}:
A \subset V(X_i), 0 < |A| \leq |V(X_i)|/2 \right \}.$$
Here, $V(X_i)$ denotes the vertex set of $X_i$, and
for a subset $A$ of $V(X_i)$, $\partial A$ denotes the
set of edges with one endpoint in $A$ and one not in $A$.

Then $G$ has {\sl Property $(\tau)$} with
respect to $\{ G_i \}$ if $\inf_i h(X(G/G_i;S)) > 0$.
This turns out not to depend on the choice of finite
generating set $S$.

When a group $G$ has Property $(\tau)$ with respect to an infinite
collection of finite index subgroups $\{ G_i \}$, there are lots
of nice consequences. For example, the resulting Schreier
coset graphs have applications in theoretical computer
science and coding theory. But when a group does not
have Property $(\tau)$, there are other useful conclusions
one can often draw, which we now briefly describe.
This is particularly the case in low-dimensional topology,
where the following important conjecture remains unresolved.

\noindent {\bf Conjecture.} [10] (Lubotzky-Sarnak) {\sl Let $G$
be the fundamental group of a finite-volume hyperbolic
3-manifold. Then $G$ does not have Property $(\tau)$ with
respect to some collection $\{G_i \}$ of finite index
subgroups.}

To appreciate the significance of this conjecture,
note the following theorem, which appears as Corollary 7.4 in [8].

\noindent {\bf Theorem.} [8] {\sl The Lubotzky-Sarnak
conjecture and the Geometrisation Conjecture together
imply that any arithmetic lattice in ${\rm PSL}(2, {\Bbb C})$
is large.}

Thus, it is important to develop new methods for showing
that a group does not have Property $(\tau)$. The following
result, which is the main theorem in this section, may be a useful tool.

\noindent {\bf Theorem 4.2.} {\sl Let $K$ be a finite
connected 2-complex with fundamental group $G$. Let $\{ G_i \}$
be a collection of finite index subgroups of $G$, and let $\{ K_i \}$
be the corresponding covering spaces of $K$. Suppose
that there is a prime $p$ and, for each $i$, there is a non-trivial class $\alpha_i$
in $H^1(K_i; {\Bbb F}_p)$, such that ${\rm relsize}(\alpha_i) \rightarrow 0$
as $i \rightarrow \infty$.
Let $\tilde G_i$ be the kernel of the homomorphism
$G_i \rightarrow {\Bbb Z}/p{\Bbb Z}$ induced by $\alpha_i$.
Then $G$ does not have Property $(\tau)$ with respect to
$\{ \tilde G_i \}$.}

In the proof of this theorem, we will need the following
construction, which will also be important later in the paper.
Let $K$ be a finite connected 2-complex with some 0-cell as a basepoint $b$.
Let $c$ be a cocycle on $K$ representing a non-trivial
element of $H^1(K; {\Bbb F}_p)$ and let $(\tilde K, \tilde b)$ be a (based) covering
space of $(K,b)$. Suppose that $\pi_1(\tilde K, \tilde b)$ lies in the
kernel of the homomorphism $\pi_1(K,b) \rightarrow {\Bbb Z}/p{\Bbb Z}$
determined by $c$. Then one can define,
for any 0-cell $v$ of $\tilde K$, its {\sl $c$-value} $c(v)$,
which is an integer mod $p$. It is defined 
as follows. Pick a path $\beta$ from $\tilde b$ to $v$ in the 1-skeleton of $\tilde K$, 
project it to a path in $K$ and define $c(v)$ to be the evaluation of $c$ on
this path. To see that this is independent of the choice of $\beta$,
let $\beta'$ be any other path in $\tilde K$ from $\tilde b$ to $v$ in the
1-skeleton of $\tilde K$.
Then $\beta'.\beta^{-1}$ is a closed loop in $\tilde K$. This projects
to a closed loop in $K$. By our hypothesis on the covering space
$\tilde K$, the evaluation of $c$ on any such closed
loop is trivial. This immediately implies that $c(v)$ is
indeed well-defined. 

An example is useful. Let $K$ be the wedge of 3 circles, labelled $x$, $y$ and $z$. Let
$c$ be the mod 2 cocycle supported on the
$x$ labelled edge. Let $\tilde K$ be the covering space corresponding
to the second term of the derived 2-series of $K$.
This is shown in Figure 1.
(Note that the dotted edges join up with each other.)
The $c$-value of its vertices is shown.

\vskip 18pt
\epsfbox{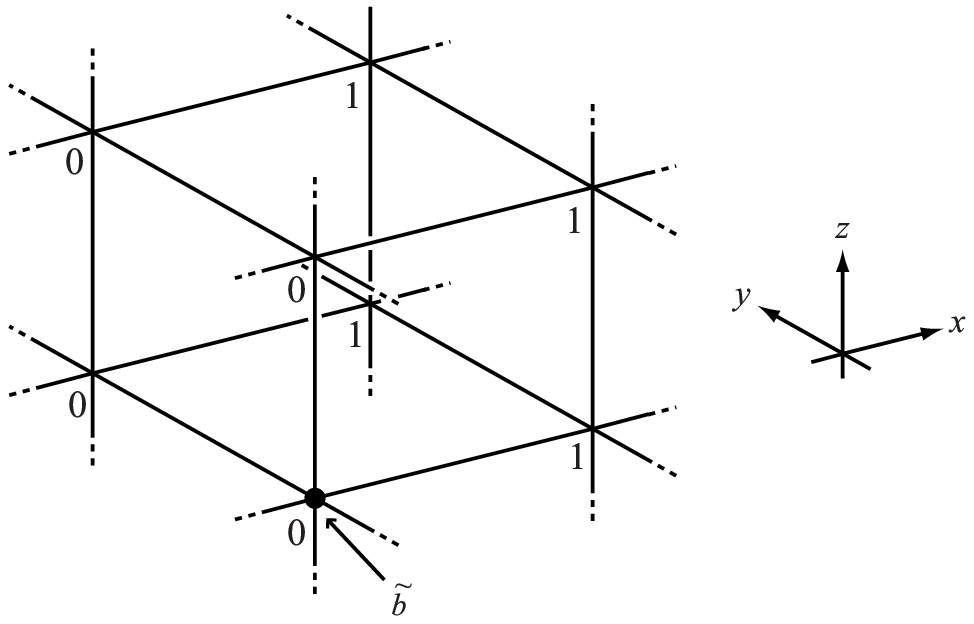}
\vskip 18pt
\centerline{Figure 1.}

\noindent {\bf Lemma 4.3.} {\sl Let $K$ be a finite connected 2-complex,
and let $\alpha$ be a non-trivial element of $H^1(K; {\Bbb F}_p)$.
Let $\tilde K$ be a finite-sheeted covering space of $K$,
such that $\pi_1(\tilde K)$ lies in the kernel of the
homomorphism $\pi_1(K) \rightarrow {\Bbb Z}/p{\Bbb Z}$ determined
by $\alpha$.
Let $\tilde X$ be the 1-skeleton of $\tilde K$. Let $V(K)$ and $E(K)$
be the 0-cells and 1-cells of $K$ respectively. Then
$$h(\tilde X) \leq {|E(K)| \over |V(K)|/p} \hbox{relsize}(\alpha).$$
}

\noindent {\sl Proof.} Let $c$ be a cocycle on $K$ representing $\alpha$
and for which $|{\rm supp}(c)|$ is minimal.
Let $A$ be the vertices in $\tilde K$ with
zero $c$-value. We claim that $|A| = |V(\tilde K)|/p$.
Let $\ell$ be any loop in $K$ based at the basepoint $b$ such that
$c(\ell) = 1$. Then $[\ell]$ represents an element
of the covering group $\pi_1(K)/\pi_1(\tilde K)$, which 
increases the $c$-value of
every vertex by 1 modulo $p$. Hence, the number of
vertices with given $c$-value is $|V(\tilde K)|/p$,
which proves the claim. As a consequence, $|A|
 = d |V(K)|/p$,
where $d$ is the degree of the cover $\tilde K \rightarrow K$.
Any edge in $\partial A$ must lie in the inverse image of the support of $c$.
Thus, $|\partial A| \leq d|{\rm supp}(c)| = d|E(K)| {\rm relsize}(\alpha)$.
So,
$$h(\tilde X) \leq {|\partial A| \over |A|} \leq
{d|E(K)| \over d|V(K)|/p} \hbox{relsize}(\alpha),$$
and the required formula now follows. $\square$

\vfill\eject
\noindent {\sl Proof of Theorem 4.2.} Let $\tilde X_i$ be the 1-skeleton of the covering
space of $K$ corresponding to $\tilde G_i$. By Lemma 4.3, 
$$h(\tilde X_i) \leq {|E(K_i)| \over |V(K_i)|/p} \hbox{relsize}(\alpha_i) = 
{|E(K)| \over |V(K)|/p} \hbox{relsize}(\alpha_i).$$
Since we are assuming that the relative size of $\alpha_i$ tends
to zero, and since the other terms on the right-hand side of the above
formula depend only on $K$, we deduce that $h(\tilde X_i) \rightarrow 0$.
Hence, $G$ does not have Property $(\tau)$ with respect to $\{ \tilde G_i \}$.
$\square$

\vskip 18pt
\centerline{\caps 5. Cocycles in covering spaces}
\vskip 6pt

Our aim over the next few chapters is to prove (2) $\Rightarrow$ (1)
of Theorems 1.14 and 1.15, and thereby establish that the group $G$ in these
theorems is large or $p$-large as appropriate. The proof will be topological, and so we consider
a finite connected 2-complex $K$ with fundamental group $G$. We are 
assuming that $G$ has a finite index subgroup $G_1$ with a 
rapidly descending $p$-series $G_i$. (In the proof of Theorem 1.15,
take $G_1$ to be $G$.) Let $K_i$ be the corresponding covering spaces
of $K$. Theorem 4.1 gives a criterion for establishing
that $G_1$ is $p$-large, in terms of the existence of 1-cocycles
on $K_i$ with relative size tending to zero. We therefore,
in this section, investigate how 1-cocycles on a 2-complex
can be used to construct 1-cocycles in covering spaces
(with potentially smaller relative size). If $U$ is a set of
1-cocycles on a cell complex $K$, we define its {\sl support}
${\rm supp}(U)$ to be the union of the supports of the cocycles
in $U$. Our main result is the following.

\noindent {\bf Theorem 5.1.} {\sl Let $K$ be a finite
connected 2-complex with $r$ 2-cells. Let $U$ be a collection
of cocycles on $K$ that represent linearly independent
elements of $H^1(K; {\Bbb F}_p)$. Let
$u = |U|$. Let $q \colon \tilde K \rightarrow K$
be a finite regular cover such that $\pi_1(K)/\pi_1(\tilde K)$
is an elementary abelian $p$-group with rank $n$. 
Then there is a collection $\tilde U$ of cocycles
on $\tilde K$ representing linearly independent
elements of $H^1(\tilde K;{\Bbb F}_p)$ such that
\item{1.} ${\rm supp}(\tilde U) \subset q^{-1}({\rm supp}(U))$;
\item{2.} $|\tilde U| \geq (n-u)u - r$.

}

The point behind Theorem 5.1 is that it provides
not just a lower bound on the dimension of $H^1(\tilde K; {\Bbb F}_p)$ 
but also gives information about certain cocycles on $\tilde K$
representing this cohomology.

We now embark on the proof of Theorem 5.1.
Consider the following subspaces of $H^1(K; {\Bbb F}_p)$:
\item{1.} the space spanned by the elements of
$U$;
\item{2.} the classes that have trivial evaluation
on all elements of $\pi_1(\tilde K)$.

\noindent Let $V_1$ and $V_2$ be these two subspaces.
Then the dimensions of $V_1$ and $V_2$ are $u$ and
$n$ respectively. 

Pick a complementary subspace for $V_1 \cap V_2$ in $V_2$,
and let $C$ be a set of cocycles on $K$ that
represents a basis for this subspace. Note that
$|C| \geq n - u$.
Note also that, by construction, $C \cup U$ forms a linearly
independent set of elements in $H^1(K; {\Bbb F}_p)$.

For each $c_1 \in U$ and $c_2 \in C$, we will show how to 
construct a cochain on $\tilde K$,
which we denote $c_1 \wedge c_2$. 
These cochains will play a vital role in the proof
of Theorem 5.1. 

Pick an orientation on each of the 1-cells
of $K$. This pulls back to give an orientation on each 1-cell $e$
of $\tilde K$. Let $i(e)$ denote the initial vertex of $e$.

Let $\tilde c_1$ be the inverse image of $c_1$ in
$\tilde K$. This is a cocycle on $\tilde K$.
Since each 1-cell $e$ is oriented, $\tilde c_1(e)$
is a well-defined element of ${\Bbb F}_p$.

Fix a basepoint $b$ in the 0-skeleton of $K$,
and let $\tilde b$ be a basepoint for $\tilde K$
in the inverse image of $b$.
Recall from Section 4 that each vertex $v$ of $\tilde K$
then has a well-defined $c_2$-value, denoted by $c_2(v)$.

We now define $c_1 \wedge c_2$.
Since the edges of $\tilde K$ are oriented, it suffices to 
assign an integer $(c_1 \wedge c_2)(e)$
modulo $p$ to each edge $e$. We define this to be
$$(c_1 \wedge c_2)(e) = \tilde c_1(e) \cdot c_2(i(e)),$$
where the product is multiplication in ${\Bbb F}_p$.

Note that ${\rm supp}(c_1 \wedge c_2) \subset 
{\rm supp}(\tilde c_1) \subset q^{-1}({\rm supp}(U))$.

It may be helpful to consider the case where
$K$ is the wedge of 3 circles labelled $x$, $y$ and $z$,
where $p = 2$ and where $\pi_1(\tilde K)$ is
the second term in the derived 2-series for $\pi_1(K)$.
Let $c_1$ and $c_2$ be the cochains supported on
the $x$-labelled and $y$-labelled edges of $K$, respectively.
Then the following
is a diagram of $\tilde K$, and the edges in 
the support of $c_1 \wedge c_2$ are
shown in bold.

\vskip 18pt
\epsfbox{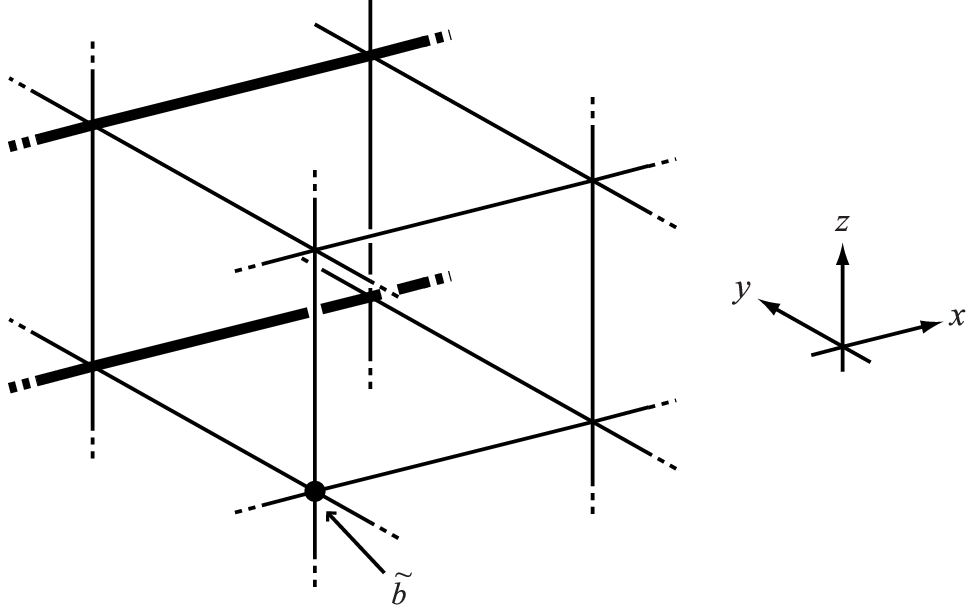}
\vskip 18pt
\centerline{Figure 2.}

We denote by $\langle U \wedge C \rangle$
the space of cochains on $\tilde K$
spanned by elements $c_1 \wedge c_2$, where $c_1 \in U$
and $c_2 \in C$. Let $Z^1(\tilde K)$ denote
the space of 1-cocycles on $\tilde K$ with mod $p$ coefficients.
We will establish the following.

\noindent {\bf Claim 5.2.} {\sl The dimension of
$\langle U \wedge C \rangle$ is $|U| \, |C|$, which is at least $u(n-u)$.}

\noindent {\bf Claim 5.3.} {\sl The subspace $Z^1(\tilde K) \cap 
\langle U \wedge C \rangle$ of cocycles 
in $\langle U \wedge C \rangle$ has codimension
at most $r$ (the number of 2-cells of $K$).}

\noindent {\bf Claim 5.4.} {\sl The map from $Z^1(\tilde K) \cap 
\langle U \wedge C \rangle$ into $H^1(\tilde K; {\Bbb F}_p)$
is an injection.}

Thus, setting $\tilde U$ to be a basis
for $Z^1(\tilde K) \cap \langle U \wedge C \rangle$
will establish Theorem 5.1.

\noindent {\bf Lemma 5.5.} {\sl Let $g$ and $h$ be
loops in $K$ based at the same point. Let $[g,h]$ denote
any lift of $ghg^{-1}h^{-1}$ to $\tilde K$. Then 
$$(c_1 \wedge c_2)([g,h]) = 
c_1(h) c_2(g) - c_1 (g) c_2(h),$$
where equality is in ${\Bbb F}_p$.}

\noindent {\sl Proof.} Each letter in the word $g$
corresponds to an edge $e$, say, in the first part of the 
loop $[g,h]$.
The inverse of this letter appears in $g^{-1}$,
where the loop runs over the edge $e'$.
Between the vertices $i(e)$ and $i(e')$ is a word $w$
conjugate to $h$. We claim that $c_2(i(e')) - c_2(i(e)) = c_2(h)$.
Pick a path from the basepoint $\tilde b$ of $\tilde K$
to $i(e)$. Then $c_2(i(e))$ is the evaluation under $c_2$
of the projection of this path to $K$. If we extend
this path using the word $w$, we obtain a path from
$\tilde b$ to $i(e')$. Thus, the difference between
$c_2(i(e'))$ and $c_2(i(e))$ is the evaluation
of $c_2$ on the projection of $w$. The parts of 
$g$ and $g^{-1}$ in $w$ project to the same edges
in $K$, but with reverse orientations. Hence,
$c_2(i(e')) - c_2(i(e))$ equals $c_2(h)$,
as required. Therefore, the evaluation of the loop
$[g,h]$ along the edges in $g$ and $g^{-1}$ is
$- c_1(g) c_2(h)$ in total. Similarly, along the edges
in $h$ and $h^{-1}$, it is $c_1(h) c_2(g)$.
So, the total evaluation is $c_1(h) c_2(g) - c_1 (g)
c_2(h)$, as required. $\square$

Let $C^1(\tilde K)$ and $B^1(\tilde K)$ be the space
of 1-cochains on $\tilde K$, with mod $p$ coefficients,
and the subspace of coboundaries.

\noindent {\bf Lemma 5.6.} {\sl The cochains $\{ c_1 \wedge c_2 :c_1 \in U, c_2 \in C\}$
map to linearly independent elements in $C^1(\tilde K)/B^1(\tilde K)$.}

\noindent {\sl Proof.}
Since $U \cup C$ forms a linearly independent set
of classes in $H^1(K; {\Bbb F}_p)$, there are
loops $\ell_i$ in $K$, based at the basepoint of $K$,
where $i \in U \cup C$,
such that, for all $c \in U \cup C$,
$$c(\ell_i) = \cases{1 & if $i = c$ \cr
0 & otherwise.}$$
Let $i \in U$ and $j \in C$. Then, by Lemma 5.5,
for any $c_1 \in U$ and $c_2 \in C$,
$$(c_1 \wedge c_2)([\ell_j, \ell_i])
= c_1(\ell_i) c_2(\ell_j) - c_1 (\ell_j)
c_2(\ell_i) =
\cases{1 &if $i=c_1$ and $j=c_2$;\cr
0 &otherwise.}$$
Since every element of $B^1(\tilde K)$ has trivial evaluation
on any loop in $\tilde K$, we deduce the lemma. $\square$

Lemma 5.6 implies Claim 5.2. It also implies that the restriction of
the quotient homomorphism $C^1(\tilde K) \rightarrow C^1(\tilde K)/B^1(\tilde K)$ 
to $\langle U \wedge C \rangle$ is an injection.
Thus, it is an injection on any subspace of $\langle U \wedge C \rangle$.
This gives Claim 5.4. We now verify Claim 5.3.

\noindent {\bf Lemma 5.7.} {\sl Let $\ell$ and
$\ell'$ be the boundary loops of two 2-cells of
$\tilde K$ that differ by a covering transformation
of $\tilde K$. Then,
$$(c_1 \wedge c_2)(\ell') = (c_1 \wedge c_2)(\ell).$$}

\noindent {\sl Proof.} Let $g$ be a path in the 1-skeleton of $\tilde K$ from
the basepoint of $\ell$ to the basepoint of $\ell'$. Thus, the loop $g\ell g^{-1}$
runs from the basepoint of $\ell$ to the basepoint of $\ell'$, then goes
around $\ell'$ and then returns to the basepoint of $\ell$.
Since $g$ and $\ell$ project to loops in $K$ based at the
same point, Lemma 5.5 gives that
$$(c_1 \wedge c_2)([g,\ell]) = 
c_1(\ell) c_2(g) - c_1 (g)c_2(\ell).$$ 
This is zero because $\ell$ is the boundary of a 2-cell
and so has zero evaluation under the cocycles
$c_1$ and $c_2$. So,
$$(c_1 \wedge c_2)(\ell') = (c_1 \wedge c_2)([g,\ell])
+ (c_1 \wedge c_2)(\ell) = (c_1 \wedge c_2)(\ell).$$
$\square$

The cocycles in $\langle U \wedge C \rangle$ are precisely those cochains
in $\langle U \wedge C \rangle$ that have zero evaluation on the boundary of any
2-cell in $\tilde K$. But Lemma 5.7 states
that if two 2-cells differ by a covering transformation,
then they have the same evaluation. Thus, one
need only check the evaluation of the boundary of
just one 2-cell in each orbit of the covering action.
There are precisely $r$ such orbits, where
$r$ is the number of 2-cells in $K$. Thus, the
codimension of $Z^1(\tilde K) \cap \langle U \wedge C \rangle$
in $\langle U \wedge C \rangle$ is at most $r$.

This proves Claim 5.3 and hence Theorem 5.1. $\square$

\noindent {\bf Remark 5.8.} Although Theorem 5.1 suffices for the purposes
of this paper, it is possible to strengthen it a little.
One can in fact find a set $\tilde U$
satisfying the requirements of Theorem 5.1, but with the
stronger inequality
$$|\tilde U| \geq un - {d (d+1) \over 2} - r,$$
where $d = {\rm dim}(V_1 \cap V_2)$. This is proved
as follows. Pick a basis for $V_1 + V_2$ so
that it contains a basis for $V_1 \cap V_2$, a basis
for $V_1$ and a basis for $V_2$. Pick a total order
on the basis for $V_1 \cap V_2$. Then consider all
cochains $c_1 \wedge c_2$, where $c_1$ lies in the
basis for $V_1$, $c_2$ lies in the basis for $V_2$,
and $c_1 < c_2$ if $c_1$ and $c_2$ both lie in $V_1 \cap V_2$.
The number of such cochains is $un - d(d+1)/2$.
It is possible to prove the corresponding versions of 5.2,
5.3 and 5.4 for these cochains. Hence, the required
inequality follows.

\noindent {\bf Remark 5.9.} The cochains $c_1 \wedge c_2$ we have considered in this section
are, in fact, special cases of a much more general construction.
In [7], a more general class of cochain was used to provide new
lower bounds on the homology growth and subgroup growth of certain
groups. These more general cochains had a certain integer $\ell$, known as their level,
assigned to them. The cochains $c_1 \wedge c_2$ are those with
level one.

\vskip 18pt
\centerline{\caps 6. The subspace reduction theorem}
\vskip 6pt

Let $E$ be a finite set, and let ${\Bbb F}_p^E$ be
the vector space over ${\Bbb F}_p$ consisting of
functions $E \rightarrow {\Bbb F}_p$.
The {\sl support} of an element $\phi$ of
${\Bbb F}_p^E$ is
$${\rm supp}(\phi) = \{e \in E : \phi(e) \not = 0 \}.$$
The {\sl support} of a subspace $W$ of ${\Bbb F}_p^E$
is 
$${\rm supp}(W) = \bigcup_{\phi \in W} {\rm supp}(\phi).$$

The main example we will consider is where $E$ is the
set of 1-cells in a finite 2-complex $K$ (with some given orientations).
Then ${\Bbb F}_p^E$ is just $C^1(K)$, the space of
1-cochains on $K$. Recall that our goal is to find
cocycles representing non-trivial elements of
$H^1(K; {\Bbb F}_p)$ and with small relative size.
The following result, which is the main theorem of
this chapter, will be the tool we use.

\noindent {\bf Theorem 6.1.} {\sl Let $V$ be a subspace
of ${\Bbb F}_p^E$ with dimension $v$, and let $w$ be a positive integer
strictly less $v$. Then, $V$ contains
a subspace $W$ with dimension $w$ such
that
$$|{\rm supp}(W)| \leq {p^v - p^{v-w} \over p^v - 1} |{\rm supp}(V)|
\leq {p^{w+1} - p \over p^{w+1} - 1}
|{\rm supp}(V)|.$$
}

In our case, $V$ will be a subspace of $C^1(K)$
spanned by $v$ cocycles, representing linearly
independent elements of $H^1(K; {\Bbb F}_p)$.
We will use Theorem 6.1 to pass to a set of
$w$ cocycles (where $w$ is a fixed integer less than $v$)
spanning a subspace $W$ with support which is
smaller than the support of $V$ by a definite factor,
independent of $v$.

We now embark on the proof of Theorem 6.1. The
following lemma gives a formula relating the
support of a subspace to the support of each
of its elements.

\noindent {\bf Lemma 6.2.} {\sl For a non-zero subspace $W$ of ${\Bbb F}_p^E$,
$$|{\rm supp}(W)| = {1 \over (p-1)p^{{\rm dim}(W) - 1}} \sum_{\phi \in W}
|{\rm supp}(\phi)|.$$}

\noindent {\sl Proof.} Focus on an element $e \in E$ in the
support of $W$. Then, there is a $\psi$ in $W$
such that $\psi(e) \not= 0$. Decompose $W$ as a direct
sum $\langle \psi \rangle \oplus W'$. Then 
we may express $W$ as a union of translates of
$W'$, as follows:
$$W = \bigcup_{i=0}^{p-1} (i \psi + W').$$
Now, for any element $\phi' \in W'$,
$i \psi(e) + \phi'(e)= 0$ for exactly one value
of $i$ between $0$ and $p-1$. Denote the indicator function of 
an element $\phi$ in ${\Bbb F}_p^E$ by $I_\phi \colon E \rightarrow \{ 0,1 \}$.
This is defined as follows:
$$I_\phi(e) = \cases{ 0 & if $\phi(e) = 0$; \cr 1 & otherwise.}$$
Then,
$$\left( \sum_{\phi \in W} I_\phi \right) (e) = 
\left( \sum_{\phi' \in W'} \sum_{i=0}^{p-1}
I_{\phi' + i \psi}\right) (e) = (p-1) p^{{\rm dim}(W) - 1}.$$
Summing this over all $e$ in the support of $W$ gives
$$\sum_{\phi \in W} |{\rm supp}(\phi)| = (p-1)p^{{\rm dim}(W) - 1}
|{\rm supp}(W)|,$$
as required. $\square$

\noindent {\bf Theorem 6.3.} {\sl Let $V$ be a non-zero subspace of
${\Bbb F}_p^E$ with dimension $v$. Then there is a codimension
one subspace $W$ of $V$ such that
$$|{\rm supp}(W)| \leq {p^v - p \over p^v - 1} |{\rm supp}(V)|.$$}

\noindent {\sl Proof.} Note that the theorem holds trivially if $v=1$.
We therefore assume $v \geq 2$.
There are $(p^v - 1)/(p-1)$ codimension one
subspaces $W$ of $V$. Summing the formula of Lemma 6.2 over
each of these gives:
$$\sum_W |{\rm supp}(W)| = {1 \over (p-1)p^{v-2}} \sum_W\sum_{\phi \in W}
|{\rm supp}(\phi)|.$$
The number of times a non-zero element $\phi$ of 
$V$ appears in the sum $\sum_W\sum_{\phi \in W}$ is
independent of the element $\phi$. Since there
are $(p^v - 1)/(p-1)$ codimension one subspaces $W$,
each containing $p^{v-1} - 1$ non-zero elements,
and there are $p^v - 1$ non-zero elements of $V$,
the number of times a non-zero element $\phi$ of 
$V$ appears in the sum $\sum_W\sum_{\phi \in W}$
is therefore
$${(p^v - 1) \over (p-1)}{ (p^{v-1} - 1) \over (p^v-1)} = 
{p^{v-1} - 1 \over p-1}.$$
Hence, 
$$\sum_W |{\rm supp}(W)| = 
{1 \over (p-1)p^{v - 2}} {(p^{v-1} - 1) \over (p-1)}
\sum_{\phi \in V} |{\rm supp}(\phi)|.$$
By Lemma 6.2,
$$\sum_{\phi \in V} |{\rm supp}(\phi)| = (p-1)p^{v-1} |{\rm supp}(V)|.$$
Hence, 
$$\sum_W |{\rm supp}(W)| = {(p^v - p) \over (p-1)} |{\rm supp}(V)|.$$
The average, over all $W$, of $|{\rm supp}(W)|$ is therefore
$${p^v - p \over p^v  -1}|{\rm supp}(V)|.$$
Hence, there is a codimension
one subspace $W$ with support at most this size. $\square$

\noindent {\sl Proof of Theorem 6.1.} 
We prove the first inequality by induction on 
the codimension $v-w$.
When this quantity is 1, this is Theorem 6.3. For the inductive step,
suppose that we have a subspace $W'$ of $V$ with dimension
$w+1$ such that
$$|{\rm supp}(W')| \leq {p^v - p^{v-w-1} \over p^v - 1}|{\rm supp}(V)|.$$
By Theorem 6.3, $W'$ has a subspace $W$ with dimension $w$ such that
$$\eqalign{
|{\rm supp}(W)| & \leq {p^{w+1} - p \over p^{w+1} - 1}{\rm supp}(W') \cr
&\leq \left({p^{w+1} - p \over p^{w+1} - 1}\right)\left( {p^v - p^{v-w-1} \over p^v - 1}\right)|{\rm supp}(V)| \cr
&={p^v - p^{v-w} \over p^v - 1} |{\rm supp}(V)|,}$$
as required. The second inequality is trivial, because
$${p^v - p^{v-w} \over p^v -1} \bigg/ {p^{w+1} - p \over p^{w+1} - 1} 
= {p^{v-w-1}(p^{w+1} -1) \over p^v  -1} = {p^v - p^{v-w-1} \over p^v -1} 
\leq 1.$$
$\square$

It is instructive to consider the case $w=1$ in Theorem 6.1.
This states that in any subspace $V$ of ${\Bbb F}_p^E$ with dimension $v>1$, 
there is an element with at most 
$${p^v- p^{v-1} \over p^v - 1}|E| \leq {p^2-p \over p^2-1} |E|$$
non-zero co-ordinates.
This is a theorem in the theory of error-correcting
codes, known as the `Plotkin bound' [14]. For, a {\sl linear code} is just a subspace
of ${\Bbb F}_p^E$, and the {\sl Hamming distance} of such a code is
the minimal number of non-zero co-ordinates in any non-zero element
of the subspace. Thus, Theorem 6.1 can be viewed as a generalisation
of the Plotkin bound, giving information not just about elements
of $V$ but whole subspaces. It is probably well-known to experts on
error-correcting codes.

\vskip 18pt
\centerline{\caps 7. Proof of Theorems 1.14 and 1.15}
\vskip 6pt

One direction of Theorems 1.14 and 1.15 is easy: the implication
$(1) \Rightarrow (2)$. The proof is as follows. Suppose that $\phi \colon G_1 \rightarrow F$
is a surjective homomorphism from a finite index
subgroup of $G$ onto a non-abelian free group $F$.
For the proof of Theorem 1.15, assume in addition
that $G_1$ is normal in $G$ and has index a power of $p$.
Let $\{F_i\}$ be the derived $p$-series of $F$, and
let $G_i = \phi^{-1}(F_i)$. Since $\{ F_i \}$
is an abelian $p$-series for $F$ with rapid descent, $\{G_i\}$ 
is therefore an abelian $p$-series with rapid descent, 
as required.

The difficult part of Theorems 1.14 and 1.15 is the implication
$(2) \Rightarrow (1)$. So, suppose that
some finite index subgroup $G_1$ of $G$
has an abelian $p$-series $\{ G_i \}$
with rapid descent. In the proof of 1.15, take $G_1$ to be $G$.
We will show that $G_1$ is $p$-large, which will establish the theorems.
Since $d_p(G_i/G_{i+1})\leq d_p(G_i)$, the rapid descent
of $\{ G_i \}$ implies that it has linear growth of mod $p$
homology. 

Let $K$ be a connected finite 2-complex with fundamental
group $G$. Let $K_i$ be the finite-sheeted covering
space corresponding to the subgroup $G_i$.
Recall from Section 4 the definition of the relative
size of an element of $H^1(K_i; {\Bbb F}_p)$, and
the following result.

\noindent {\bf Theorem 4.1.} {\sl Let $K$ be a finite connected
2-complex, and let $\{ K_i \rightarrow K\}$ be a
collection of finite-sheeted covering spaces.
Suppose that $\{ \pi_1(K_i) \}$ has linear
growth of mod $p$ homology for some prime $p$.
Then one of the following must hold:
\item{(i)} $\pi_1(K_i)$ is $p$-large for infinitely many $i$, or
\item{(ii)} there is some $\epsilon > 0$ such that
the relative size of any non-trivial class in $H^1(K_i; {\Bbb F}_p)$
is at least $\epsilon$, for all $i$.

}

Thus, our plan is to prove that (ii) of Theorem 5.1 does not hold,
and therefore deduce that $G_1$ is $p$-large. We will keep track of
a set $U_i$ of cellular 1-dimensional cocycles on $K_i$ that represent
linearly independent elements of $H^1(K_i; {\Bbb F}_p)$.
The cardinality $|U_i|$ will be
some fixed positive integer $u$ independent of $i$.
(The precise size of $u$ will depend on data from
the group $G$ and the series $\{ G_i \}$.)
Our aim is to ensure that
$${ |{\rm supp}(U_i)| \over \hbox{Number of 1-cells of } K_i}
\rightarrow 0.\eqno{(\dag)}$$
In particular, the relative size of any
element of $U_i$ tends to zero, which means that
(ii) does not hold.

We establish $(\dag)$ using the following method.
Let $q_i \colon K_{i+1} \rightarrow K_{i}$ be
the covering map. We will find a set of cocycles
$U_{i+1}^+$ on $K_{i+1}$ representing linearly 
independent elements of $H^1(K_{i+1}; {\Bbb F}_p)$, 
with the following properties:
\item{I.} ${\rm supp}(U_{i+1}^+) \subset q_{i}^{-1}({\rm supp}(U_i))$;
\item{II.} $|U_{i+1}^+| > u$.

\noindent Note that the inequality in (II) is
strict. 

Let $E$ denote the set of 1-cells of $K_{i+1}$ with given orientation.
Then $C^1(K_{i+1})$ is 
isomorphic to ${\Bbb F}_p^E$, the vector space of
functions $E \rightarrow {\Bbb F}_p$. 
Let $V$ be the subspace of $C^1(K_{i+1})$ spanned
by $U_{i+1}^+$, and let $w = u$. We apply Theorem 6.1 to $V$.

\noindent {\bf Theorem 6.1.} {\sl Let $V$ be a subspace
of ${\Bbb F}_p^E$ with dimension $v$, and let $w$ be a positive integer
strictly less $v$. Then, $V$ contains
a subspace $W$ with dimension $w$ such
that
$$|{\rm supp}(W)| \leq {p^v - p^{v-w} \over p^v - 1} |{\rm supp}(V)|
\leq {p^{w+1} - p \over p^{w+1} - 1}
|{\rm supp}(V)|.$$
}

Let $U_{i+1}$ be a basis for the subspace $W$ given
by Theorem 6.1. 
Note that the factor $(p^{w+1} - p)/( p^{w+1} - 1)$
is strictly less than 1, and is dependent only on
the fixed integers $u$ and $p$. Now, $|q_i^{-1}({\rm supp}(U_i))|$
is obtained from $|{\rm supp}(U_i)|$ by scaling by the
degree of the cover $q_i$. The same relation holds between
the number of 1-cells in $K_{i+1}$ and the number of 1-cells in
$K_i$. Thus, $(\dag)$ follows.

The key, then, is to construct the cocycles $U_{i+1}^+$
with properties (I) and (II).
For this, we use Theorem 5.1 (setting
$K = K_i$, $U = U_i$ and $\tilde K = K_{i+1}$)

\noindent {\bf Theorem 5.1.} {\sl Let $K$ be a finite
connected 2-complex with $r$ 2-cells. Let $U$ be a collection
of cocycles on $K$ that represent linearly independent
elements of $H^1(K; {\Bbb F}_p)$. Let
$u = |U|$. Let $q \colon \tilde K \rightarrow K$
be a finite regular cover such that $\pi_1(K)/\pi_1(\tilde K)$
is an elementary abelian $p$-group with rank $n$. 
Then there is a \break collection $\tilde U$ of cocycles
on $\tilde K$ representing linearly independent
elements of \break $H^1(\tilde K; {\Bbb F}_p)$ such that
\item{1.} ${\rm supp}(\tilde U) \subset q^{-1}({\rm supp}(U))$;
\item{2.} $|\tilde U| \geq (n-u)u - r$.

}

We define $U_{i+1}^+$ to be the set $\tilde U$ provided by
this theorem. Condition (1) of the theorem is just
(I) above. We need to ensure that
Condition (II) holds. Thus, we require
$$(n-u)u - r > u.$$
We will ensure that this holds by using the hypothesis
that $G_i$ has rapid descent and by a suitable choice of $u$.

Let $\langle X | R \rangle$ be a finite presentation for $G$.
Let $\lambda$ be 
$$\liminf_i {d_p(G_i/G_{i+1}) - 2 \over [G:G_i]}.$$
Since $\{ G_i \}$ is rapidly descending, $\lambda$ is 
positive. Let $u$ be 
$$\left \lceil {4|R| \over \lambda} \right \rceil,$$
which is a positive integer.
Now pick a sufficiently large integer $j$, such that, for all $i \geq j$,
$${d_p(G_i/G_{i+1}) - 2 \over [G:G_i]} > {\lambda \over 2}$$
and 
$$\lambda[G:G_i] > 4u.$$
Hence,
$$\lambda[G:G_i]/2 - u > \lambda[G:G_i]/4.$$
Now, $$n = d_p(G_i/G_{i+1}) \geq \lambda[G:G_i]/2 +2.$$
So,
$$(n-u)u - r \geq
(\lambda[G:G_i]/2 + 2-u)u - r
\geq \lambda [G:G_i]u/4 + 2u - r
\geq 2u,$$
since $$\lambda u \geq 4|R| = 4r/[G:G_i].$$
This proves (2) $\Rightarrow$ (1) of Theorems 1.14 and 1.15. $\square$

\vskip 18pt
\centerline{\caps References}
\vskip 6pt

\item{1.} {\caps G. Baumslag}, {\sl Groups with the same 
lower central sequence as a relatively free group. I. The groups.} 
Trans. Amer. Math. Soc. 129 (1967) 308--321. 

\item{2.} {\caps G. Baumslag}, {\sl Groups with the same lower central sequence as 
a relatively free group. II. Properties.} Trans. Amer. Math. Soc. 142 (1969) 507--538.

\item{3.} {\caps D. Cooper, D. Long, A. Reid,} 
{\sl Essential closed surfaces in bounded $3$-manifolds,} 
J. Amer. Math. Soc. 10 (1997) 553--563.

\item{4.} {\caps M. Lackenby}, {\sl A characterisation of
large finitely presented groups}, 
J. Algebra. 287 (2005) 458--473.

\item{5.} {\caps M. Lackenby}, {\sl Covering spaces of 3-orbifolds},
Duke Math. J. 136 (2007) 181--203. 

\item{6.} {\caps M. Lackenby}, {\sl Large groups, Property $(\tau)$ 
and the homology growth of subgroups}, Preprint.

\item{7.} {\caps M. Lackenby}, {\sl New lower bounds on subgroup growth
and homology growth}, Preprint.

\item{8.} {\caps M. Lackenby, D. Long, A. Reid}, {\sl Covering spaces
of arithmetic 3-orbifolds}, Preprint.

\item{9.} {\caps M. Lackenby,} {\sl Pro-$p$ groups and low-dimensional topology}, 
In preparation.

\item{10.} {\caps A. Lubotzky}, {\sl Eigenvalues of the Laplacian, the 
First Betti Number and the Congruence Subgroup Problem},
Ann. Math. 144 (1996) 441--452.

\item{11.} {\caps A. Lubotzky, D. Segal},
{\sl Subgroup growth}. Progress in Mathematics, 212. 
Birkh\"auser Verlag (2003)

\item{12.} {\caps W. L\"uck}, {\sl Approximating $L^2$-invariants by their
finite-dimensional analogues}, G.A.F.A. 4 (1994) 455--481.

\item{13.} {\caps W. L\"uck}, {\sl $L^2$-invariants of regular coverings of compact manifolds 
and CW-complexes,} Handbook of geometric topology, Elsevier (2002)

\item{14.} {\caps M. Plotkin}, {\sl Binary codes with specified minimum distance},
IRE Transations on Information Theory 6 (1960) 445--450.

\item{15.} {\caps Y. Shalom}, {\sl
Rigidity, unitary representations of semisimple groups, 
and fundamental groups of manifolds with rank one transformation group.}
Ann. Math.  152 (2000) 113--182.

\vskip 12pt
\+ Mathematical Institute, University of Oxford, \cr
\+ 24-29 St Giles', Oxford OX1 3LB, United Kingdom. \cr

\end